\documentclass{scrartcl}

\usepackage{tikz-cd}
\usepackage{adjustbox}
\usepackage{amsmath,amsfonts,amssymb,amsthm, stmaryrd}
\usepackage{hyperref}
\usepackage{esint}
\usepackage{float}
\usepackage{doi}
\usepackage{braket}
\usepackage{quiver}
\usepackage{mama}

\theoremstyle{definition}

\theoremstyle{theorem}

\newcommand{\FibOptic}{\cat{FibOptic}}

\newcommand{\DLens}{\cat{DLens}}
\newcommand{\Lens}{\cat{Lens}}

\newcommand{\DMark}{\cat{DMark}}
\newcommand{\Optic}{\cat{Optic}}

\newcommand{\Para}{\cat{Para}}
\newcommand{\Copara}{\cat{Copara}}
\newcommand{\CoPara}{\Copara}

\newcommand{\DOptic}{\cat{DOptic}}
\newcommand{\Mat}{\cat{Mat}}
\newcommand{\dom}{\mathrm{dom}}

\newcommand{\mor}{\mathrm{mor}}
\newcommand{\deloop}{\twocat B}

\newcommand{\Ca}{{\mathcal C}}
\newcommand{\Da}{{\mathcal D}}
\newcommand{\Ma}{{\mathcal M}}
\newcommand{\Ia}{{\mathcal I}}

\newcommand{\action}{\bullet}
\newcommand{\laction}{\rhd}
\newcommand{\lactionn}{\blacktriangleright}
\newcommand{\raction}{\lhd}

\newcommand{\forwards}{\mathsf{forwards}}
\newcommand{\backwards}{\mathsf{backwards}}
\newcommand{\view}{\mathsf{view}}
\newcommand{\update}{\mathsf{update}}
\newcommand{\copyfunc}{\mathsf{copy}}

\newcommand{\Fr}{\twocat{F}\cat{r}}

\newcommand{\groth}{\int}

\usepackage{calc}
\newlength{\hght}
\newlength{\dpth}

\title{Fibre optics \\ (extended abstract)\footnote{This is an extended abstract for a forthcoming paper of the same title.}}
\author{Dylan Braithwaite \and Matteo Capucci \and Bruno Gavranovi\'c \and Jules Hedges \and Eigil Fjeldgren Rischel}
\date{}

\begin{document}

\maketitle

\section{Introduction}
\label{section:introduction}
Lenses \cite{hedges_lenses_philosophers}, optics \cite{CategoriesOfOptics,pickering_gibbons_wu_profunctor_optics,ProfunctorOptics} and dependent lenses \cite{abbott-etal-categories-containers,spivak_poly} (or equivalently morphisms of containers, or equivalently natural transformations of polynomial functors) are all widely used in applied category theory as models of bidirectional processes.
From the definition of lenses over a finite product category, optics weaken the required structure to actions of monoidal categories, and dependent lenses make use of the additional property of finite completeness (or, in case of polynomials, even local cartesian closure).
This has caused a split in the applied category theory literature between those using optics (for example \cite{bolt_hedges_zahn_bayesian_open_games,towards_foundations_categorical_cybernetics,escrows_optics}) and those using dependent lenses (for example \cite{spivak_poly,cts,smithe_compositional_active_inference_1}).
The goal of this paper is to unify optics with dependent lenses, by finding a definition of \emph{fibre optics} admitting both as special cases.
Informally, we are completing the following `square', where the arrows mean `generalises to':
\[\begin{tikzcd}[sep=3ex]
	{\text{lenses}} &&& {\text{optics}} \\
	\\
	{\text{dep. lenses}} &&& {?}
	\arrow[from=1-1, to=3-1]
	\arrow[from=1-1, to=1-4]
	\arrow[dashed, from=1-4, to=3-4]
	\arrow[dashed, from=3-1, to=3-4]
\end{tikzcd}\]

We refer the reader to \cite{ProfunctorOptics} for background and definitions on optics, and \cite{abbott-etal-categories-containers} on dependent lenses (called \emph{morphisms of containers} in that paper).

\section{Indexed optics}
\label{section:ommatidia}
One characterisation of the category of dependent lenses over $\mathbf{Set}$ is that it is the coproduct completion of lenses. (This should not be confused with the \emph{free} coproduct completion, for which the embedding does not preserve coproducts).
Coproducts are quite useful for the applications we have in mind, namely categorical cybernetics (\cite{other_paper,towards_foundations_categorical_cybernetics}), as they are the semantics of an `external choice' operation (demonic choice) on bidirectional processes.

Motivated by this, we would like to characterise the coproduct completion of categories of mixed optics.
Suppose $\Ma$ is a distributive monoidal category\footnotemark\ acting on categories $\Ca$, $\Da$, with both actions preserving coproducts in $\Ma$.
\footnotetext{A distributive monoidal category is a monoidal category with all coproducts, for which coproducts distribute over the monoidal product. For example, the Kleisli category of any commutative monad on $\Sets$ is distributive monoidal.}
An object of the resulting category is a set-indexed family of pairs of objects of $\Ca$ and $\Da$.
The set of morphisms $\left( I, \binom{X_i}{X'_i} \right) \to \left( J, \binom{Y_j}{Y'_j} \right)$ is given by the coend
\[
	\int^{M \in [I \times J, \Ma]}
	\prod_{i \in I} \Ca \left( X_i, \sum_{j \in J} M_i^j \action Y_j) \right)
	\times
	\Da \left(\sum_{j \in J} Y'_j \action M_i^j, X'_i \right)
\]
Identity morphisms are given by taking $M$ to be an identity matrix, that is, $M_i^j = I$ (the monoidal unit of $\Ma$) when $i = j$, and $0$ (the initial object of $\Ma$) otherwise. Composition is given by matrix multiplication: $(MN)_i^k = \sum_{j \in J} M_i^j \otimes N_j^k$.

This is closely related to Milewski's \emph{ommatidia} (named after the compound eyes of insects), which are less flexible but allow us to save considerable work because they arise as an ordinary category of optics \cite{milewski_polylens}.

If the actions of $\Ma$ on $\Ca$ and $\Da$ additionally have right adjoints (analogous to monoidal closed categories) then $\Ca$ and $\Da$ become enriched in $\Ma$ with the action as a tensor \cite[§2.2]{gray1980closed}.
Under these conditions, we can prove that the coproduct completion of optics is equivalent to the category of \emph{polynomial $\Ma$-functors $\Da \to \Ca$}, defined as coproducts of tensors of $\Ma$-representables: $F(Y) = \sum_{i \in I} [Y, X'_i] \bullet X_i$.

This definition is satisfactory when every object of the base category is a coproduct of copies of the monoidal unit. An example of this is the category of sets and finite support Markov kernels. In general, objects of the coproduct completion are ``piecewise trivial bundles'': bundles $X' \to \sum_{i \in I} X_i$ which are trivial when restricted to every $X_i$. This motivates the general definition in the next section, which expands the applicability of this approach beyond piecewise trivial bundles.

\section{Fibre Optics}
\label{section:indexed-optics}
Let $\Ia$ be a category with all finite limits, considered as a cartesian monoidal category.
Let $\Ma, \Ca, \Da \to \Ia$ be Beck--Chevalley monoidal bifibrations \cite{shulman_framed_bicategories}, and let $\Ma$ act on $\Ca$ and $\Da$ in a way which is compatible with the bifibrations (we call this a \emph{fibraction}).
In this case, the monoidal structure on $\Ma$ induces a monoidal structure on each fibre $\Ma_{I}$, and the actions of $\Ma$ on $\Ca$ and $\Da$ induce an action of each $\Ma_{I}$ on $\Ca_{I}$ and $\Da_{I}$.
Moreover, this fiberwise structure is respected by the pullback functors $f^*$, which are strong monoidal and preserve the action (the pushforward functors $f^!$ do not in general have either property, although they inherit an oplax monoidal structure by virtue of being left adjoint to a lax monoidal functor \cite{kelly1974doctrinal}).

Over this data we define a category $\FibOptic_\Ma^\Ia(\Ca,\Da)$ whose objects are triples $\left( I, \binom{X}{X'} \right)$ where $I \in \Ia$, $X \in \Ca_I$ and $X' \in \Da_I$.
The set of morphisms $\left( I, \binom{X}{X'} \right) \to \left( J, \binom{Y}{Y'} \right)$ is given by the coend
\begin{equation}
\label{eq:fibre-optics}
	\int^{M \in \Ma_{I \times J}}
	\Ca_I(C, \pi_I^!(M\action \pi_J^* C'))
	\times
	\Da_I(\pi_I^!(\pi_J^* D' \action M), D)
\end{equation}
where $\pi_I$ and $\pi_J$ are the projections out of $I \times J$.
This specialises to the previously defined set-indexed optics in the case where the fibre $\Ca_I$ is the category of functors $[I, \Ca]$.

As an example, let $\Meas$ be the category of measurable spaces and measurable functions. Let $\DMark$ be the category whose objects are measurable functions $A \to X$, where a morphism from $p_A : A \to X$ to $p_B : B \to Y$ consists of a Markov kernel $k : A \to B$ and a measurable function $f : X \to Y$ such that $\mathbb{P}[p_{B}(b) = f(p_{A}(a)) \mid b \sim k(a)] = 1$.

There is an obvious functor $\DMark \to \Meas$ which just remembers the codomain, and one can show that this is a Beck--Chevalley monoidal bifibration. $\DMark_{X}$ is a category of ``Markov kernels indexed over $X$''.
The ordinary product defines a monoidal structure on $\DMark$ which makes this a monoidal bifibration.
Applying the above construction for $\cat{I} = \Meas$, $\cat{C} = \cat{D} = \cat{M} = \DMark$ (acting on itself by cartesian product), one finds a category of ``general stochastic indexed lenses''.

\section{A proposal for a general definition of dependent optics}
Once we figured out the story of indexed optics and fibre optics, we started investigating whether they would fit in a more conceptual framework for `dependent optics', i.e. a general mathematical theory of optics with `linear' dependent types.
In fact, fibre optics work well when the actions involved are multiplicative\footnote{An informal term which roughly means `product-like', and that we could formalize crudely as `coproducts distributes over actions'.}, but can we set the stage for more exotic constructions?

We advance a proposal motivated by the examples we expound immediately after, in which we use it to reconstruct lenses, mixed optics and fibre optics alike.
The proposal is firmly grounded in the idea that optics are naturally decomposed in their forward and backward parts, and that residuals play an integral role in linking them.

Informally, our proposal is to define dependent optics as simply as pullbacks of bicategories:
\[\begin{tikzcd}
	\DOptic & {\text{bicat. of forward parts}} \\
	{\text{bicat. of backward parts}} & {\text{cat. of residuals}}
	\arrow[from=1-2, to=2-2]
	\arrow[""{name=0, anchor=center, inner sep=0}, from=2-1, to=2-2]
	\arrow["\backwards"', from=1-1, to=2-1]
	\arrow["\forwards", from=1-1, to=1-2]
	\arrow["\lrcorner"{anchor=center, pos=0.125}, draw=none, from=1-1, to=0]
\end{tikzcd}\]

We believe that, as stated, this definition might be too general.
Therefore, at the end of this abstract we propose a refinement, where we ask the functors projecting down residuals (the legs of the pullback) to be 2-(op)fibrations embodying some higher categorical versions of actegories.
Details are still in the workings.

\subsection{Mixed optics as dependent optics}
For a monoidal category $\Ma$ we write $\deloop\Ma$ for its delooping, which is a 1-object bicategory whose composition and coherence morphisms are borrowed from $\Ma$'s \cite[Example 2.1.19]{johnson2021}.
For a left $\Ma$-actegory $(\Da, \laction)$, $\Para_\laction(\Da)$ is the bicategory whose objects are those of $\Da$ and whose morphisms are $\Ma$-parametrised morphisms $(M : \Ma, f : M \laction X \to Y)$ \cite[Definition 2]{towards_foundations_categorical_cybernetics}.
Dually, given a left $\Ma$-actegory $(\Ca, \lactionn)$, $\Copara (\Ca)$ has morphisms $(M \in \Ma, f : X \to M \lactionn Y)$.

Most importantly, $\Para_\laction(\Da)^\co$ is naturally 2-opfibred over $\deloop\Ma$, by projecting parameters, and $\Copara_\lactionn(\Ca)$ is similarly naturally 2-fibred over $(\deloop\Ma)^\op$.

We call $\Optic_{\lactionn,\laction}$ the 2-category arising from the pullback:
\[\begin{tikzcd}
	{\Optic_{\lactionn,\laction}} & {\Para_\laction(\Da)^\coop} \\
	{\Copara_\lactionn(\Ca)} & {(\deloop\Ma)^\op}
	\arrow[""{name=0, anchor=center, inner sep=0}, "C"', from=2-1, to=2-2]
	\arrow["P", from=1-2, to=2-2]
	\arrow[dashed, from=1-1, to=2-1]
	\arrow[dashed, from=1-1, to=1-2]
	\arrow["\lrcorner"{anchor=center, pos=0.125}, draw=none, from=1-1, to=0]
\end{tikzcd}\]
In particular, one gets the usual 1-category of optics (as defined in \cite{categoricalupdate2019}) by locally quotienting such a pullback by taking connected components (informally, by `changing base along $\pi_0 : \Cat \to \Sets$').

Without this quotient, the bicategory of optics we obtain remembers explicitly the choices of residuals made for an optics, and slidings remain as explicit 2-morphisms instead of disappearing under the usual coend.
Under some conditions to ensure strictness, this corresponds to replacing the coends in $\Sets$ in the classical definition with oplax coends \cite[Definition 7.1.5]{loregian-coend-cofriend} in $\Cat$.

\subsection{Lenses and dependent lenses as dependent optics}
\label{section:cosmic-tesseract}
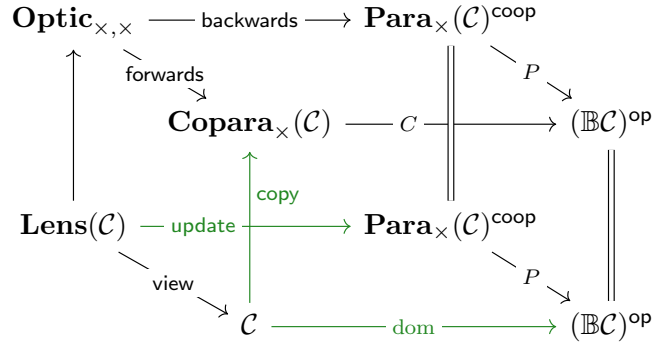
\begin{figure}
\[\begin{tikzcd}[sep=2ex, column sep=0.5ex]
	{\Optic_{\times,\times}} && {\Para_\times(\Ca)^\coop} \\
	\\
	& {\CoPara_\times(\Ca)} && {(\deloop \Ca)^\op} \\
	\\
	{\Lens(\Ca)} && {\Para_\times(\Ca)^\coop} \\
	\\
	& \Ca && {(\deloop \Ca)^\op}
	\arrow["\dom"{description}, color={rgb,255:red,34;green,135;blue,34}, from=7-2, to=7-4]
	\arrow["P"{description}, from=5-3, to=7-4]
	\arrow[Rightarrow, no head, from=7-4, to=3-4]
	\arrow["\copyfunc"'{pos=0.7}, color={rgb,255:red,34;green,135;blue,34}, from=7-2, to=3-2]
	\arrow["\view"{description}, from=5-1, to=7-2]
	\arrow[from=5-1, to=1-1]
	\arrow["\update"{description, pos=0.3}, color={rgb,255:red,34;green,135;blue,34}, from=5-1, to=5-3]
	\arrow["\forwards"{description}, from=1-1, to=3-2]
	\arrow["\backwards"{description}, from=1-1, to=1-3]
	\arrow["P"{description}, from=1-3, to=3-4]
	\arrow["C"{description, pos=0.3}, from=3-2, to=3-4]
	\arrow[Rightarrow, no head, from=5-3, to=1-3]
\end{tikzcd}\]
	\caption{The `Cosmic' Cube}
	\label{fig:cosmic_cube}
\end{figure}

Contrary to optics, lenses (and dependent lenses) come with a specific and fixed choice of residual, which is given by the domain of the view part, which gets `copied' down to the update.
Usually one proves optics for cartesian actions collapse to lenses by using coend calculus reductions on each hom-set.
The definition of optics as pullbacks provides a new perspective by showing how lenses arise when we construct optics with a fixed choice of residual (Figure~\ref{fig:cosmic_cube}).

In the cube, $\Ca$ is a cartesian monoidal category.
The bottom face shows a construction of lenses as dependent optics, by presenting them as pullback of a functor describing the coparametrisaiton of their view part and a functor describing the parametrisation of their update part.
Indeed, the $\dom$ functor is an oplax functor sending every candidate view morphism to its domain (the fixed residual of a lens).
The oplax $\copyfunc$ functor instead sends a a morphism $f:A \to B$ to its `canonically coparametrised' version $\Delta_A \comp (1_A \times f) : A \to A \times B$.
The top face is instead the construction of mixed optics expounded above, and is a diagram of bicategories.

Once we locally quotient all the bicategories, the cube becomes strictly commutative and provides a new way to see how optics for cartesian actions collapse to lenses.

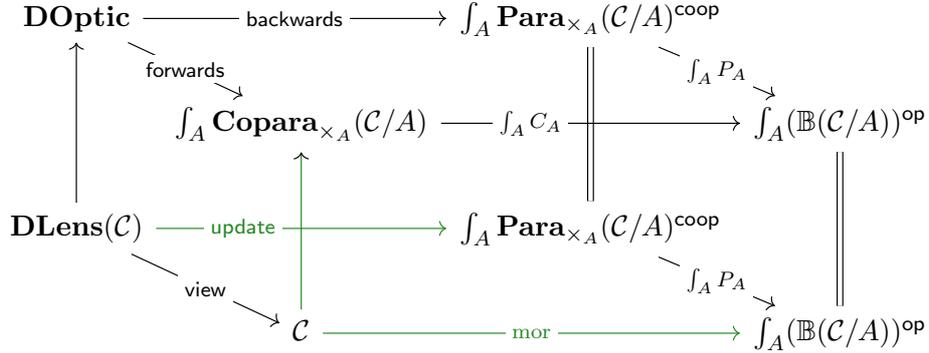
\begin{figure}[H]
\[\begin{tikzcd}[sep=2ex, column sep=0.5ex]
	\DOptic && {\groth_A \Para_{\times_A}(\Ca/A)^\coop} \\
	\\
	& {\groth_A \CoPara_{\times_A}(\Ca/A)} && {\groth_A (\deloop(\Ca/A))^\op} \\
	\\
	{\DLens(\Ca)} && {\groth_A \Para_{\times_A}(\Ca/A)^\coop} \\
	\\
	& \Ca && {\groth_A (\deloop(\Ca/A))^\op}
	\arrow["\view"{description}, from=5-1, to=7-2]
	\arrow["\mor"{description}, color={rgb,255:red,34;green,135;blue,34}, from=7-2, to=7-4]
	\arrow["{\groth_A P_A}"{description}, from=5-3, to=7-4]
	\arrow["\update"{description, pos=0.3}, color={rgb,255:red,34;green,135;blue,34}, from=5-1, to=5-3]
	\arrow[from=5-1, to=1-1]
	\arrow[Rightarrow, no head, from=5-3, to=1-3]
	\arrow["\backwards"{description}, from=1-1, to=1-3]
	\arrow[Rightarrow, no head, from=7-4, to=3-4]
	\arrow[draw={rgb,255:red,34;green,135;blue,34}, from=7-2, to=3-2]
	\arrow["{\int_A C_A}"{description, pos=0.3}, from=3-2, to=3-4]
	\arrow["{\int_A P_A}"{description}, from=1-3, to=3-4]
	\arrow["\forwards"{description}, from=1-1, to=3-2]
\end{tikzcd}\]
	\caption{The `Cosmic' Dependent Cube}
	\label{fig:cosmic_dependent_cube}
\end{figure}

We also conjecture that for a finitely complete category $\Ca$, dependent lenses in $\Ca$ fit in a similar diagram (Figure~\ref{fig:cosmic_dependent_cube}).
The main idea is to perform fiberwise para and copara constructions on the fibers of the indexed cartesian monoidal category $\Ca/-$.

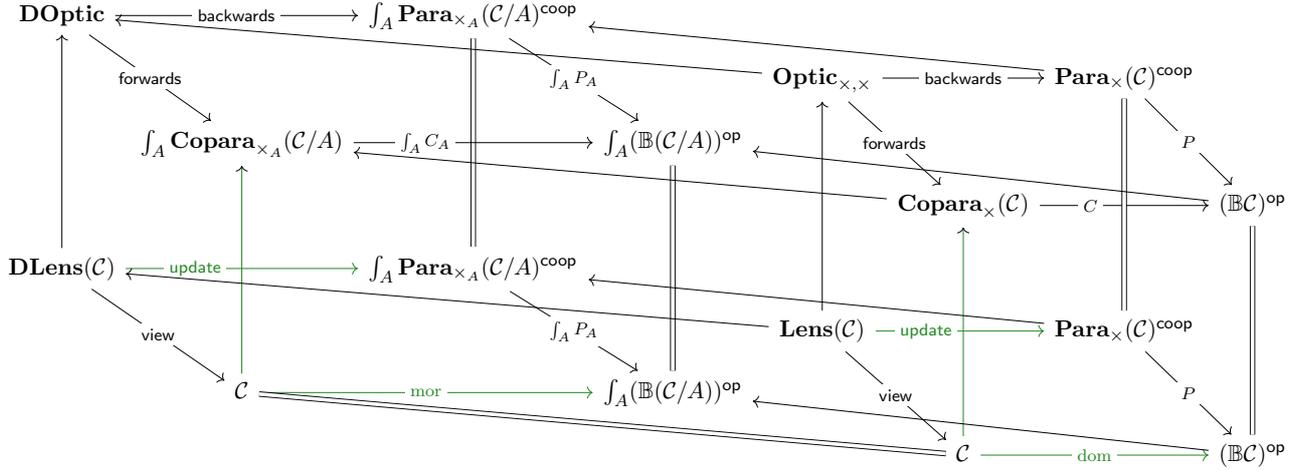
\begin{figure}
	\[\adjustbox{scale=0.8,center}{\begin{tikzcd}[sep=2ex, column sep=0.5ex]
	\DOptic && {\groth_A \Para_{\times_A}(\Ca/A)^\coop} \\
	&&&&& {\Optic_{\times,\times}} && {\Para_\times(\Ca)^\coop} \\
	& {\groth_A \CoPara_{\times_A}(\Ca/A)} && {\groth_A (\deloop(\Ca/A))^\op} \\
	&&&&&& {\CoPara_\times(\Ca)} && {(\deloop \Ca)^\op} \\
	{\DLens(\Ca)} && {\groth_A \Para_{\times_A}(\Ca/A)^\coop} \\
	&&&&& {\Lens(\Ca)} && {\Para_\times(\Ca)^\coop} \\
	& \Ca && {\groth_A (\deloop(\Ca/A))^\op} \\
	&&&&&& \Ca && {(\deloop \Ca)^\op}
	\arrow["\dom"{description}, color={rgb,255:red,34;green,135;blue,34}, from=8-7, to=8-9]
	\arrow["P"{description}, from=6-8, to=8-9]
	\arrow[Rightarrow, no head, from=8-9, to=4-9]
	\arrow[draw={rgb,255:red,34;green,135;blue,34}, from=8-7, to=4-7]
	\arrow["\view"{description}, from=6-6, to=8-7]
	\arrow[from=6-6, to=2-6]
	\arrow["\update"{description, pos=0.3}, color={rgb,255:red,34;green,135;blue,34}, from=6-6, to=6-8]
	\arrow["\forwards"{description}, from=2-6, to=4-7]
	\arrow["\backwards"{description}, from=2-6, to=2-8]
	\arrow["P"{description}, from=2-8, to=4-9]
	\arrow["C"{description, pos=0.3}, from=4-7, to=4-9]
	\arrow["\view"{description}, from=5-1, to=7-2]
	\arrow["\mor"{description}, color={rgb,255:red,34;green,135;blue,34}, from=7-2, to=7-4]
	\arrow["{\groth_A P_A}"{description}, from=5-3, to=7-4]
	\arrow["\update"{description, pos=0.3}, color={rgb,255:red,34;green,135;blue,34}, from=5-1, to=5-3]
	\arrow[from=5-1, to=1-1]
	\arrow[Rightarrow, no head, from=5-3, to=1-3]
	\arrow["\backwards"{description}, from=1-1, to=1-3]
	\arrow[Rightarrow, no head, from=7-4, to=3-4]
	\arrow[draw={rgb,255:red,34;green,135;blue,34}, from=7-2, to=3-2]
	\arrow[Rightarrow, no head, from=7-2, to=8-7]
	\arrow[Rightarrow, no head, from=6-8, to=2-8]
	\arrow["{\int_A C_A}"{description, pos=0.3}, from=3-2, to=3-4]
	\arrow["{\int_A P_A}"{description}, from=1-3, to=3-4]
	\arrow["\forwards"{description}, from=1-1, to=3-2]
	\arrow[from=6-6, to=5-1]
	\arrow[from=8-9, to=7-4]
	\arrow[from=6-8, to=5-3]
	\arrow[from=4-7, to=3-2]
	\arrow[from=4-9, to=3-4]
	\arrow[from=2-6, to=1-1]
	\arrow[from=2-8, to=1-3]
\end{tikzcd}}\]
	\caption{The Cosmic Tesseract}
	\label{fig:cosmic_tesseract}
\end{figure}

Finally, the two cubes can be tied together by embedding simple lenses in dependent lenses, thereby producing a tesseract (Figure~\ref{fig:cosmic_tesseract}).

\subsection{Fibre optics as dependent optics}
Next we show that our explicit definition of fibre optics is also an instance of dependent optics.

We start by noticing there is a bicategory $\Mat_\Ia (\Ma)$ whose objects are objects of $\Ia$, whose 1-cells $I \to J$ are objects $M: \Ma_{I \times J}$ with 2-morphisms given by morphisms between them in $\Ma_{I \times J}$.
This is the category $\Fr(\Ma)$ defined in \cite[Section 14]{shulman_framed_bicategories}, where it is presented as a framed bicategory.
It can be seen as a generalisation of $\deloop\Ma$, which one gets when $\Ia$ is the terminal category.

Then we can define an indexed category $\cat{MatMult}_\Ca: \Mat_\Ia(\Ma)^\op \to \Cat$ sending an object $I$ to the fibre $\Ca_I$ of $\Ca$, and morphisms $M : I \to J$ to the `matrix-vector multiplication' functor
\[
	\pi_I^! (M \action \pi_J^* (-)): \Ca_J \longto \Ca_I.
\]
From this we define a fibred version of $\Para$ and $\CoPara$:
\begin{align*}
	&\CoPara_\Ma(\Ca) = \int_{I \in \Mat_\Ia(\Ma)}\cat{MatMult}_\Ca(I)\\
	&\Para_\Ma(\Da) = \left( \int_{I \in \Mat_\Ia(\Ma)} \cat{MatMult}_\Da(I)^\op\right)^\op
\end{align*}
An object of $\Copara_\Ma(\Ca)$ is a pair $(I, C)$ where $I \in \Ia$ and $C \in \Ca_I$, and a morphism $(I, C) \to (J, C')$ is a morphism $C \to \pi_I^!(M \action \pi_J^* C')$ in $\Ca_I$ where $M \in \Ma_{I \times J}$.
Similarly a morphism $(I, D) \to (J, D')$ in $\Para_\Ma(\Da)$ is a morphism $\pi_I^!(\pi_J^* D \action M) \to D'$ in $\C_I$.
It's therefore clear~\eqref{eq:fibre-optics} can be reproduced analogously to the way we reproduced mixed optics:
\[\begin{tikzcd}
	\FibOptic_\Ma^\Ia(\Ca,\Da) & {\CoPara_\Ma(\Ca)} \\
	{\Para_\Ma(\Da)^\op} & \Mat_\Ia(\Ma)^\op
	\arrow[from=1-2, to=2-2]
	\arrow[from=2-1, to=2-2]
	\arrow[from=1-1, to=2-1]
	\arrow[from=1-1, to=1-2]
	\arrow["\lrcorner"{anchor=center, pos=0.125}, draw=none, from=1-1, to=2-2]
\end{tikzcd}\]
The category $\Mat_\Ia(\Ma)$ generalises from the action of a one-object bicategory in the case of optics, to the action of a `many-objects' bicategory.
This action is similar to the action of a bicategory defined in \cite{bakovic2009simplicial} once instanced for a framed bicategory \cite{shulman_framed_bicategories}.
With such an action comes two \emph{action bicategories}, $\Ca \raction \Ma$ and $\Ma \laction \Da$ which generalise $\CoPara$ and $\Para$, respectively, and are likewise 2-fibred and 2-opfibred over $\Ma$.
This would yield a pullback diagram
\[\begin{tikzcd}
	{\DOptic_{\Ma}(\Ca, \Da)} & {\Ca \raction \Ma} \\
	{(\Ma^\op \laction \Da)^\op} & \Ma
	\arrow[from=1-2, to=2-2]
	\arrow[from=2-1, to=2-2]
	\arrow[from=1-1, to=2-1]
	\arrow[from=1-1, to=1-2]
	\arrow["\lrcorner"{anchor=center, pos=0.125}, draw=none, from=1-1, to=2-2]
\end{tikzcd}\]
which we think of the most general instance of dependent optics still grounded on the `indexed forward/backward' intuition.
We are still working out whether any instance of this, apart from the ones we already know, yields interesting constructions (we speculate $\Delta$-lenses might arise in this way).

\bibliographystyle{alpha}
\bibliography{../references}

\end{document}